\documentclass[12pt]{article}
\usepackage{amssymb,amsmath}
\usepackage{hyperref}
\usepackage{graphicx}
\usepackage{color}

\begin{document}

\title{\LARGE\bf A reformulated series expansion of the arctangent function}

\author{
\normalsize\bf S. M. Abrarov\footnote{\scriptsize{Dept. Earth and Space Science and Engineering, York University, Toronto, Canada, M3J 1P3.}}\, and B. M. Quine$^{*}$\footnote{\scriptsize{Dept. Physics and Astronomy, York University, Toronto, Canada, M3J 1P3.}}}

\date{January 15, 2017}
\maketitle

\begin{abstract}
In our recent publication we obtained a series expansion of the arctangent function involving complex numbers. In this work we show that this formula can also be expressed as a real rational function.
\vspace{0.25cm}
\vspace{0.25cm}
\end{abstract}

\vspace{0.25cm}
\section{Derivation}

\subsection{Derivatives of the arctangent function}

\vspace{0.25cm}
Consider the following function
$$
\frac{-i}{{{\left( x+i \right)}^{m}}},
$$
where $x$ is a real variable and $m$ is a non-negative integer. Since
\small
$$
\frac{-i}{{{\left( x+i \right)}^{1}}} = -\frac{1}{1+{{x}^{2}}}-i\frac{x}{1+{{x}^{2}}}
$$
it is not difficult to show by induction that
\newpage
\small
$$
\begin{aligned}
\frac{-i}{{{\left( x+i \right)}^{m}}}&=\frac{1}{{{\left( x+i \right)}^{m-1}}}\frac{-i}{{{\left( x+i \right)}^{1}}} =\frac{1}{{{\left( x+i \right)}^{m-1}}}\left( -\frac{1}{1+{{x}^{2}}}-i\frac{x}{1+{{x}^{2}}} \right) \\ 
& =i^{1}\left( \frac{-i}{{{\left( x+i \right)}^{m-1}}} \right)\left( -\frac{1}{1+{{x}^{2}}}-i\frac{x}{1+{{x}^{2}}} \right)^1 \\
& ={{i}^{2}}\left( \frac{-i}{{{\left( x+i \right)}^{m-2}}} \right){{\left( -\frac{1}{1+{{x}^{2}}}-i\frac{x}{1+{{x}^{2}}} \right)}^{2}} =\cdots \\ 
& ={{i}^{n}}\left( \frac{-i}{{{\left( x+i \right)}^{m-n}}} \right){{\left( -\frac{1}{1+{{x}^{2}}}-i\frac{x}{1+{{x}^{2}}} \right)}^{n}}=\cdots \\
& =-{{i}^{m+1}}{{\left( -\frac{1}{1+{{x}^{2}}}-i\frac{x}{1+{{x}^{2}}} \right)}^{m}}, \qquad\qquad\qquad m \geq n.  
\end{aligned}
$$
\normalsize
Applying the binomial formula to this identity
\footnotesize
\[
-{{i}^{m+1}}{{\left( -\frac{1}{1+{{x}^{2}}}-i\frac{x}{1+{{x}^{2}}} \right)}^{m}}=-{{i}^{m+1}}\sum\limits_{n=0}^{m}{{{\left( -\frac{1}{1+{{x}^{2}}} \right)}^{m-n}}{{\left( -\frac{ix}{1+{{x}^{2}}} \right)}^{n}}\left( 
\begin{aligned}
  & m \\ 
 & n \\ 
\end{aligned} \right)}
\]
\normalsize
and then separating its right side into the real and imaginary parts, after some trivial rearrangements we obtain
\footnotesize
\begin{equation}\label{eq_1}
\frac{-i}{{{\left( x+i \right)}^{m}}}=\sum\limits_{n=1}^{m}{\frac{{{\left( -1 \right)}^{n}}{{x}^{m-\left( 2n-1 \right)}}}{{{\left( 1+{{x}^{2}} \right)}^{m}}}\left( \begin{matrix}
  m \\ 
  2n-1 \\ 
\end{matrix} \right)+}\,\,i\frac{{{\left( -1 \right)}^{n}}{{x}^{m-2\left( n-1 \right)}}}{{{\left( 1+{{x}^{2}} \right)}^{m}}}\left( \begin{matrix}
  m \\ 
  2\left( n-1 \right) \\ 
\end{matrix} \right).
\end{equation}
\normalsize
Taking into consideration that
$$
\operatorname{Re}\left[ \frac{-i}{{{\left( x+i \right)}^{m}}} \right]=-\operatorname{Re}\left[ \frac{-i}{{{\left( x-i \right)}^{m}}} \right]
$$
and
$$
\operatorname{Im}\left[ \frac{-i}{{{\left( x+i \right)}^{m}}} \right]=\operatorname{Im}\left[ \frac{-i}{{{\left( x-i \right)}^{m}}} \right]
$$
we can also find
\footnotesize
\begin{equation}\label{eq_2}
-\frac{-i}{{{\left( x-i \right)}^{m}}}=\sum\limits_{n=1}^{m}{\frac{{{\left( -1 \right)}^{n}}{{x}^{m-\left( 2n-1 \right)}}}{{{\left( 1+{{x}^{2}} \right)}^{m}}}\left( \begin{matrix}
  m \\ 
  2n-1 \\ 
\end{matrix} \right)-}\,\,i\frac{{{\left( -1 \right)}^{n}}{{x}^{m-2\left( n-1 \right)}}}{{{\left( 1+{{x}^{2}} \right)}^{m}}}\left( \begin{matrix}
  m \\ 
  2\left( n-1 \right) \\ 
\end{matrix} \right).
\end{equation}
\normalsize
Consequently, taking sum of equations \eqref{eq_1} and \eqref{eq_2} provides
\begin{equation}\label{eq_3}
\begin{aligned}
\left( \frac{-i}{{{\left( x+i \right)}^{m}}}-\frac{-i}{{{\left( x-i \right)}^{m}}} \right)&=\frac{1}{i}\left( \frac{1}{{{\left( x+i \right)}^{m}}}-\frac{1}{{{\left( x-i \right)}^{m}}} \right)\\
&=2\sum\limits_{n=1}^{m}{\frac{{{\left( -1 \right)}^{n}}{{x}^{m-\left( 2n-1 \right)}}}{{{\left( 1+{{x}^{2}} \right)}^{m}}}\left( \begin{matrix}
  m \\ 
  2n-1 \\ 
\end{matrix} \right)}.
\end{aligned}
\end{equation}

The ${{m}^{\text{th}}}$ derivative of the arctangent function can be represented as \cite{Dinghui1999, Abrarov2016a} 
\begin{equation}\label{eq_4}
\frac{{{d}^{m}}}{d\,{{x}^{m}}}\arctan \left( x \right)=\frac{{{\left( -1 \right)}^{m}}\left( m-1 \right)!}{2i}\left( \frac{1}{{{\left( x+i \right)}^{m}}}-\frac{1}{{{\left( x-i \right)}^{m}}} \right).
\end{equation}
Comparing equations \eqref{eq_3} and \eqref{eq_4} immediately yields 
\begin{equation}\label{eq_5}
\frac{{{d}^{m}}}{d\,{{x}^{m}}}\arctan \left( x \right)=\left( m-1 \right)!\sum\limits_{n=1}^{m}{\frac{{{\left( -1 \right)}^{m+n}}{{x}^{m-\left( 2n-1 \right)}}}{{{\left( 1+{{x}^{2}} \right)}^{m}}}\left( \begin{matrix}
  m \\ 
  2n-1 \\ 
\end{matrix} \right)}.
\end{equation}
Notably, in contrast to the equation \eqref{eq_4}, the series expansion \eqref{eq_5} involves no complex numbers.

It interesting to note that the equation \eqref{eq_4} is also directly related to the identity
\[
\frac{{{d^m}}}{{d{t^m}}}\arctan \left( x \right) = {\rm{sgn}}^{m - 1}\left( { - x} \right)\frac{{\left( {m - 1} \right)!}}{{{{\left( {1 + {x^2}} \right)}^{m/2}}}}\sin \left( {m\arcsin \left( {\frac{1}{{\sqrt {1 + {x^2}} }}} \right)} \right),
\]
where the signum function is defined as
\[
{\rm{sgn}}\left( x \right) = \left\{
\begin{aligned}
1, & \qquad x \ge 0\\
 - 1, & \qquad x < 0,
\end{aligned}
\right.
\]
by de Moivre's formula (see \cite{Dinghui1999} for derivation procedure). However, in the work \cite{Dinghui1999} this identity is shown without signum function that, according to Lampret \cite{Lampret2011}, is necessary to account for $x \in \mathbb{R^{-}}$ in derivatives of the acrtangent function.

\subsection{The arctangent function}

In our recent publication we have derived a series expansion of the arctangent function \cite{Abrarov2016b}
\small
\begin{equation}\label{eq_6}
\arctan \left( x \right)=i\underset{M\to \infty }{\mathop{\lim }}\,\sum\limits_{m=1}^{\left\lfloor \frac{M}{2} \right\rfloor +1}{\frac{1}{2m-1}\left( \frac{1}{{{\left( 1+2i/x \right)}^{2m-1}}}-\frac{1}{{{\left( 1-2i/x \right)}^{2m-1}}} \right)}.
\end{equation}
\normalsize
We used the notation $\left\lfloor M/2 \right\rfloor +1$ in equation \eqref{eq_6} only for chronological reason to keep consistency with its previously published variation (see \cite{Abrarov2016a} for details)
\[
\begin{aligned}
\arctan \left( x \right)=&\,\,i\underset{L\to \infty }{\mathop{\lim }}\,\sum\limits_{\ell =1}^{L}\sum\limits_{m=1}^{\left\lfloor \frac{M}{2} \right\rfloor +1}\frac{1}{2m-1} \\
&\times \left( \frac{1}{{{\left( \left( 2\ell -1 \right)+2i/x \right)}^{2m-1}}}-\frac{1}{{{\left( \left( 2\ell -1 \right)-2i/x \right)}^{2m-1}}} \right).
\end{aligned}
\]

Since at $M\to \infty $  the upper limit in summation $\left\lfloor M/2 \right\rfloor +1$ also tend to infinity, the series expansion \eqref{eq_6} can be rewritten as
\begin{equation}\label{eq_7}
\arctan \left( x \right)=i\sum\limits_{m=1}^{\infty }{\frac{1}{2m-1}\left( \frac{1}{{{\left( 1+2i/x \right)}^{2m-1}}}-\frac{1}{{{\left( 1-2i/x \right)}^{2m-1}}} \right)}.
\end{equation}

In order to exclude the complex numbers in equation \eqref{eq_7} we can apply a similar approach as we have made already for derivation of the identity \eqref{eq_5}. In particular, we note that
\begin{equation}\label{eq_8}
\frac{1}{{{\left( 1+2i/x \right)}^{m}}}=\frac{{{\left( x/2 \right)}^{m}}}{-i}\frac{-i}{{{\left( x/2+i \right)}^{m}}}.
\end{equation}
Consequently, for the function
$$
\frac{-i}{{{\left( x/2+i \right)}^{m}}}
$$
we can simply replace the variable $x\to x/2$ in the identity \eqref{eq_1} to obtain
\footnotesize
\begin{equation}\label{eq_9}
\begin{aligned}
&\frac{-i}{{{\left( x/2+i \right)}^{m}}} \\
& \hspace{1cm}=\sum\limits_{n=1}^{m}{\frac{{{\left( -1 \right)}^{n}}{{\left( \frac{x}{2} \right)}^{m-\left( 2n-1 \right)}}}{{{\left( 1+{{\left( \frac{x}{2} \right)}^{2}} \right)}^{m}}}\left( \begin{matrix}
  m \\ 
  2n-1 \\ 
\end{matrix} \right)+}\,\,i\frac{{{\left( -1 \right)}^{n}}{{\left( \frac{x}{2} \right)}^{m-2\left( n-1 \right)}}}{{{\left( 1+{{\left( \frac{x}{2} \right)}^{2}} \right)}^{m}}}\left( \begin{matrix}
  m \\ 
  2\left( n-1 \right) \\ 
\end{matrix} \right).
\end{aligned}		
\end{equation}
\normalsize
Comparing now equations \eqref{eq_8} and \eqref{eq_9} we get
\vspace{0.2cm}
\footnotesize
\begin{equation}\label{eq_10}
\hspace{-11.5cm}\frac{1}{{{\left( 1+2i/x \right)}^{m}}}
\end{equation}
\[
\qquad=i{{\left( \frac{x}{2} \right)}^{m}}\left( \sum\limits_{n=1}^{m}{\frac{{{\left( -1 \right)}^{n}}{{\left( \frac{x}{2} \right)}^{m-\left( 2n-1 \right)}}}{{{\left( 1+{{\left( \frac{x}{2} \right)}^{2}} \right)}^{m}}}\left( \begin{matrix}
  m \\ 
  2n-1 \\ 
\end{matrix} \right)+i\frac{{{\left( -1 \right)}^{n}}{{\left( \frac{x}{2} \right)}^{m-2\left( n-1 \right)}}}{{{\left( 1+{{\left( \frac{x}{2} \right)}^{2}} \right)}^{m}}}\left( \begin{matrix}
  m \\ 
  2\left( n-1 \right) \\ 
\end{matrix} \right)} \right).	
\]
\normalsize

Similarly, writing
$$
\frac{1}{{{\left( 1-2i/x \right)}^{m}}}=\frac{{{\left( x/2 \right)}^{m}}}{-i}\frac{-i}{{{\left( x/2-i \right)}^{m}}}
$$
and replacing $x\to x/2$ in equation \eqref{eq_2} results in
\footnotesize
\begin{equation}\label{eq_11}
\hspace{-11.5cm}-\frac{1}{{{\left( 1-2i/x \right)}^{m}}} \\
\end{equation}
\[
\qquad = i{{\left( \frac{x}{2} \right)}^{m}} \left( \sum\limits_{n=1}^{m}{\frac{{{\left( -1 \right)}^{n}}{{\left( \frac{x}{2} \right)}^{m-\left( 2n-1 \right)}}}{{{\left( 1+{{\left( \frac{x}{2} \right)}^{2}} \right)}^{m}}}\left( \begin{matrix}
  m \\ 
  2n-1 \\ 
\end{matrix} \right)-i\frac{{{\left( -1 \right)}^{n}}{{\left( \frac{x}{2} \right)}^{m-2\left( n-1 \right)}}}{{{\left( 1+{{\left( \frac{x}{2} \right)}^{2}} \right)}^{m}}}\left( \begin{matrix}
  m \\ 
  2\left( n-1 \right) \\ 
\end{matrix} \right)} \right).
\]
\normalsize
Therefore, taking sum of equations \eqref{eq_10} and \eqref{eq_11} leads to
$$
\frac{1}{{{\left( 1+2i/x \right)}^{m}}}-\frac{1}{{{\left( 1-2i/x \right)}^{m}}}=2i{{\left( \frac{x}{2} \right)}^{m}}\sum\limits_{n=1}^{m}{\frac{{{\left( -1 \right)}^{n}}{{\left( \frac{x}{2} \right)}^{m-\left( 2n-1 \right)}}}{{{\left( 1+{{\left( \frac{x}{2} \right)}^{2}} \right)}^{m}}}\left( \begin{matrix}
  m \\ 
  2n-1 \\ 
\end{matrix} \right)}
$$
or
\begin{equation}\label{eq_12}
\begin{aligned}
\frac{1}{m}\left( \frac{1}{{{\left( 1+2i/x \right)}^{m}}}\right. & \left.-\frac{1}{{{\left( 1-2i/x \right)}^{m}}} \right) \\
& = \frac{2i}{m}{{\left( \frac{x}{2} \right)}^{m}}\sum\limits_{n=1}^{m}\frac{{{\left( -1 \right)}^{n}}{{\left( \frac{x}{2} \right)}^{m-\left( 2n-1 \right)}}}{{{\left( 1+{{\left( \frac{x}{2} \right)}^{2}} \right)}^{m}}}\left( \begin{matrix}
  m \\ 
  2n-1 \\ 
\end{matrix} \right).
\end{aligned}
\end{equation}

Replacing in equation \eqref{eq_12} the index $m\to 2m-1$ provides
\small
\begin{equation}\label{eq_13}
\begin{aligned}
& \frac{1}{2m-1}\left( \frac{1}{{{\left( 1+2i/x \right)}^{2m-1}}} -\frac{1}{{{\left( 1-2i/x \right)}^{2m-1}}} \right) \\
& \hspace{2cm}=\frac{2i}{2m-1}{{\left( \frac{x}{2} \right)}^{2m-1}}\sum\limits_{n=1}^{2m-1}{\frac{{{\left( -1 \right)}^{n}}{{\left( \frac{x}{2} \right)}^{2m-1-\left( 2n-1 \right)}}}{{{\left( 1+{{\left( \frac{x}{2} \right)}^{2}} \right)}^{2m-1}}}\left( \begin{matrix}
  2m-1 \\ 
  2n-1 \\ 
\end{matrix} \right)}.
\end{aligned}
\end{equation}
\normalsize
Consequently, applying the equation \eqref{eq_13} to the formula \eqref{eq_7} we have
\small
\[
\begin{aligned}
& i\sum\limits_{m=1}^{\infty }\frac{1}{2m-1}\left( \frac{1}{{{\left( 1+2i/x \right)}^{2m-1}}} - \frac{1}{{{\left( 1-2i/x \right)}^{2m-1}}} \right) \\
& \hspace{1.75cm} = i\sum\limits_{m=1}^{\infty }{\frac{2i}{2m-1}{{\left( \frac{x}{2} \right)}^{2m-1}}\sum\limits_{n=1}^{2m-1}{\frac{{{\left( -1 \right)}^{n}}{{\left( \frac{x}{2} \right)}^{2m-1-\left( 2n-1 \right)}}}{{{\left( 1+{{\left( \frac{x}{2} \right)}^{2}} \right)}^{2m-1}}}\left( \begin{matrix}
  2m-1 \\ 
  2n-1 \\ 
\end{matrix} \right)}}
\end{aligned}
\]
\normalsize
or
\small
\[
\begin{aligned}
& i\sum\limits_{m=1}^{\infty }{\frac{1}{2m-1}\left( \frac{1}{{{\left( 1+2i/x \right)}^{2m-1}}}-\frac{1}{{{\left( 1-2i/x \right)}^{2m-1}}} \right)} \\ 
& \hspace{1.75cm} = -2\sum\limits_{m=1}^{\infty }{\sum\limits_{n=1}^{2m-1}{\frac{{{\left( -1 \right)}^{n}}}{\left( 2m-1 \right){{\left( 1+{{x}^{2}}/4 \right)}^{2m-1}}}{{\left( \frac{x}{2} \right)}^{2\left( 2m-n \right)-1}}\left( \begin{matrix}
  2m-1 \\ 
  2n-1 \\ 
\end{matrix} \right)}}
\end{aligned}
\]
\normalsize
or
\begin{equation}\label{eq_14}
\hspace{-10cm}\arctan \left( x \right)= \\
\end{equation}
\[
-2\sum\limits_{m=1}^{\infty }{\sum\limits_{n=1}^{2m-1}{\frac{{{\left( -1 \right)}^{n}}}{\left( 2m-1 \right){{\left( 1+{{x}^{2}}/4 \right)}^{2m-1}}}{{\left( \frac{x}{2} \right)}^{2\left( 2m-n \right)-1}}\left( \begin{matrix}
  2m-1 \\ 
  2n-1 \\ 
\end{matrix} \right)}}.
\]
As we can see, the equation \eqref{eq_14} is represented in form of a real rational function.

Although performing computation by truncating equation \eqref{eq_14} is not optimal due to double summation, its application may be more convenient in theoretical analysis. In particular, reformulation as the real rational function \eqref{eq_14} can help understand, for example, a behavior of equation \eqref{eq_7} at small $x \to 0$ in computing pi (see \cite{Abrarov2016b} for details) and estimate its error in truncation.

\section{Conclusion}
We derived an equivalent \eqref{eq_14} to series expansion \eqref{eq_7} of the arctangent function in form a real rational function. Since equation \eqref{eq_14} involves no complex numbers, its application may be more convenient in theoretical analysis.

\section*{Acknowledgments}

This work is supported by National Research Council Canada, Thoth Technology Inc. and York University. The authors wish to thank Dr. J. Guillera and Dr. M. Shao for discussions and useful information.

\bigskip


\end{document}